\documentclass[12pt,twoside]{article}
\usepackage{amsfonts}
\usepackage{amsmath}

\date{}

\begin{document}

\date{}

\centerline{\bf Journal's Title, Vol. x, 200x, no. xx, xxx - xxx}

\centerline{}

\centerline{}

\centerline {\Large{\bf On The Special Curves in}}

\centerline{}

\centerline{\Large{\bf Minkowski 4 Spacetime}}

\centerline{}

\centerline{\bf {G\"{u}l G\"{u}ner}}

\centerline{}

\centerline{Karadeniz Technical University}

\centerline{Department of Mathematics}

\centerline{Trabzon, Turkey}

\centerline{gguner@ktu.edu.tr}

\centerline{}

\centerline{\bf {F. Nejat Ekmekci}}

\centerline{}

\centerline{Ankara University}

\centerline{Department of Mathematics}

\centerline{Ankara, Turkey}

\centerline{ekmekci@science.ankara.edu.tr}

\newtheorem{Theorem}{\quad Theorem}[section]

\newtheorem{Definition}[Theorem]{\quad Definition}

\newtheorem{Corollary}[Theorem]{\quad Corollary}

\newtheorem{Lemma}[Theorem]{\quad Lemma}

\newtheorem{Example}[Theorem]{\quad Example}

\begin{abstract}
In [1], we gave a method for constructing Bertrand curves from the spherical
curves in 3 dimensional Minkowski space. In this work, we construct the
Bertrand curves corresponding to a spacelike geodesic and a null helix in
Minkowski 4 spacetime.
\end{abstract}



\section{Preliminary Notes}

In this section, we give basic notions of spacelike and null curves in
Minkowski 4-space (see $\left[ 2\right] ,\left[ 3\right] $ and $\left[ 6%
\right] $). Let $%
\mathbb{R}
^{4}=\{(x_{1},x_{2},x_{3},x_{4}):x_{1},x_{2},x_{3},x_{4}\in
\mathbb{R}
\}$ be a 4-dimensional vector space. For any vectors $%
x=(x_{1},x_{2},x_{3},x_{4}),y=(y_{1},y_{2},y_{3},y_{4})$ in $%
\mathbb{R}
^{4}$, the pseudo scalar product of $x$ and $y$ is defined to be $%
\left\langle x,y\right\rangle =-x_{1}y_{1}+x_{2}y_{2}+x_{3}y_{3}+x_{4}y_{4}$%
. We call $(%
\mathbb{R}
^{4},\left\langle ,\right\rangle )$ a Minkowski 4-space. We write $%
\mathbb{R}
_{1}^{4}$ instead of $(%
\mathbb{R}
^{4},\left\langle ,\right\rangle )$. We say that a non-zero vector $x\in
\mathbb{R}
_{1}^{4}$ is spacelike, lightlike (null) or timelike if $\left\langle
x,x\right\rangle >0$, $\left\langle x,x\right\rangle =0$ or $\left\langle
x,x\right\rangle <0$ respectively. The norm of the vector $x\in
\mathbb{R}
_{1}^{4}$ is defined by $\left\Vert x\right\Vert =\sqrt{|\left\langle
x,x\right\rangle |}$. For a vector $v\in
\mathbb{R}
_{1}^{4}$ and a real number $c$, we define a hyperplane with pseudo normal $%
v $ by $HP(v,c)=\{x\in
\mathbb{R}
_{1}^{4}:\left\langle x,v\right\rangle =c\}.$ We call $HP(v,c)$ a spacelike
hyperplane, a timelike hyperplane or a lightlike hyperplane if $v$ is
timelike, spacelike or lightlike respectively. We also define de Sitter
3-space by $S_{1}^{3}=\{x\in
\mathbb{R}
_{1}^{4}:\left\langle x,x\right\rangle =1\}.$ For any $%
x=(x_{1},x_{2},x_{3},x_{4}),y=(y_{1},y_{2},y_{3},y_{4}),z=\left(
z_{1},z_{2},z_{3},z_{4}\right) $ in $%
\mathbb{R}
_{1}^{4}$, we define a vector%
\begin{equation*}
x\wedge y\wedge z=\left\vert
\begin{array}{cccc}
-e_{1} & e_{2} & e_{3} & e_{4} \\
x_{1} & x_{2} & x_{3} & x_{4} \\
y_{1} & y_{2} & y_{3} & y_{4} \\
z_{1} & z_{2} & z_{3} & z_{4}%
\end{array}%
\right\vert
\end{equation*}%
where $(e_{1},e_{2},e_{3},e_{4})$ is the canonical basis of $%
\mathbb{R}
_{1}^{4}$. We can easily show that $\left\langle a,(x\wedge y\wedge
z)\right\rangle =\det (a,x,y,z).$

Let $\gamma :I\longrightarrow S_{1}^{3}$ be a regular curve. We say that a
regular curve $\gamma $ is spacelike, timelike or null respectively, if $%
\gamma ^{\prime }(t)$ is spacelike, timelike or null at any $t\in I$, where $%
\gamma ^{\prime }=d\gamma /dt$. Now we describe the explicit differential
geometry on spacelike and null curves in $S_{1}^{3}$.

Let $\gamma $ be a spacelike regular curve, we can reparametrise $\gamma $
by the arclength $s=s(t)$. Hence, we may assume that $\gamma (s)$ is a unit
speed curve. So we have the tangent vector $t(s)=\gamma ^{\prime }(s)$ with $%
\left\Vert t(s)\right\Vert =1$. In the case when $\left\langle t^{\prime
}(s),t^{\prime }(s)\right\rangle \neq 1$, we have a unit vector $n(s)=\dfrac{%
t^{\prime }(s)-\gamma (s)}{\left\Vert t^{\prime }(s)-\gamma (s)\right\Vert }$%
. Moreover, define $e(s)=\gamma (s)\wedge t(s)\wedge n(s)$, then we have a
pseudo orthonormal frame $\{\gamma (s),t(s),n(s),e(s)\}$ of $%
\mathbb{R}
_{1}^{4}$ along $\gamma $. By the standard arguments, we can show the
following Frenet-Serret type formulae: Under the assumption that $%
\left\langle t^{\prime }(s),t^{\prime }(s)\right\rangle \neq 1$,%
\begin{eqnarray}
\gamma ^{\prime }\left( s\right) &=&t\left( s\right)  \notag \\
t^{\prime }\left( s\right) &=&-\gamma \left( s\right) +\kappa _{g}\left(
s\right) n\left( s\right)  \notag \\
n^{\prime }\left( s\right) &=&\kappa _{g}\left( s\right) \delta \left(
\gamma \left( s\right) \right) t\left( s\right) +\tau _{g}\left( s\right)
e\left( s\right)  \notag \\
e^{\prime }\left( s\right) &=&\tau _{g}\left( s\right) n\left( s\right)
\end{eqnarray}%
where $\delta (\gamma (s))=-sign(n(s))$,
\begin{eqnarray*}
\kappa _{g}(s) &=&\left\Vert t^{\prime }(s)+\gamma (s)\right\Vert \\
\tau _{g}\left( s\right) &=&\dfrac{\delta \left( \gamma \left( s\right)
\right) }{\kappa
{{}^2}%
_{g}\left( s\right) }\det \left( \gamma \left( s\right) ,\gamma ^{\prime
}\left( s\right) ,\gamma ^{\prime \prime }\left( s\right) ,\gamma ^{\prime
\prime \prime }\left( s\right) \right)
\end{eqnarray*}

Now let $\gamma :I\longrightarrow S_{1}^{3}$ be a null curve. We will
assume, in the sequel, that the null curve we consider has no points at
which the acceleration vector is null. Hence $\left\langle \gamma ^{\prime
\prime }\left( t\right) ,\gamma ^{\prime \prime }\left( t\right)
\right\rangle $ is never zero. We say that a null curve $\gamma \left(
t\right) $ in $%
\mathbb{R}
_{1}^{4}$ is parametrized by the pseudo-arc if $\left\langle \gamma ^{\prime
\prime }\left( t\right) ,\gamma ^{\prime \prime }\left( t\right)
\right\rangle =1.$ If a null curve satisfies $\left\langle \gamma ^{\prime
\prime }\left( t\right) ,\gamma ^{\prime \prime }\left( t\right)
\right\rangle \neq 0$, then $\left\langle \gamma ^{\prime \prime }\left(
t\right) ,\gamma ^{\prime \prime }\left( t\right) \right\rangle >0$, and
\begin{equation*}
u\left( t\right) =\overset{t}{\underset{t_{0}}{\int }}\left\langle \gamma
^{\prime \prime }\left( t\right) ,\gamma ^{\prime \prime }\left( t\right)
\right\rangle ^{1/4}dt
\end{equation*}%
becomes the pseudo-arc parameter.

A null curve $\gamma \left( t\right) $ in $%
\mathbb{R}
_{1}^{4}$ with $\left\langle \gamma ^{\prime \prime }\left( t\right) ,\gamma
^{\prime \prime }\left( t\right) \right\rangle \neq 0$ is a Cartan curve if $%
\left\{ \gamma ^{\prime }\left( t\right) ,\gamma ^{\prime \prime }\left(
t\right) ,\gamma ^{\prime \prime \prime }\left( t\right) \right\} $ is
linearly independent for any $t$ . For a Cartan curve $\gamma \left(
t\right) $ in $%
\mathbb{R}
_{1}^{4}$ with pseudo-arc parameter $t$ , there exists a pseudo orthonormal
basis $\{L,N,W_{1},W_{2}\}$ such that%
\begin{eqnarray}
L &=&\gamma ^{\prime }  \notag \\
L^{\prime } &=&W_{1}  \notag \\
N^{\prime } &=&-\gamma +k_{1}W_{1}+k_{2}W_{2}  \notag \\
W_{1}^{\prime } &=&-k_{1}L-N  \notag \\
W_{2}^{\prime } &=&-k_{2}L
\end{eqnarray}%
where $\left\langle L,N\right\rangle =1,\left\langle L,W_{1}\right\rangle
=\left\langle L,W_{2}\right\rangle =\left\langle N,W_{1}\right\rangle
=\left\langle N,W_{2}\right\rangle =\left\langle W_{1},W_{2}\right\rangle
=0. $ We call $\{L,N,W_{1},W_{2}\}$ as the Cartan frame and $\left\{
k_{1},k_{2}\right\} $ as the Cartan curvatures of $\gamma $. Since the
Cartan frame is unique up to orientation, the number of the Cartan
curvatures is minimum and the Cartan curvatures are invariant under Lorentz
transformations, the set $\{L,N,W_{1},W_{2},k_{1},k_{2}\}$ corresponds to
the Frenet apparatus of a space curve. A direct computation shows that the
values of the Cartan curvatures are%
\begin{eqnarray}
k_{1} &=&\frac{1}{2a%
{{}^2}%
}\left( \left\langle \gamma ^{\prime \prime \prime },\gamma ^{\prime \prime
\prime }\right\rangle +2aa^{\prime \prime }-4\left( a^{\prime }\right)
{{}^2}%
\right)  \notag \\
k_{2} &=&-\frac{1}{a^{4}}\det (\gamma ^{\prime },\gamma ^{\prime \prime
},\gamma ^{\prime \prime \prime },\gamma ^{\left( 4\right) })
\end{eqnarray}%
\begin{Theorem}\textit{Let }$\gamma \left( t\right) $\textit{%
\ in }$%
\mathbb{R}
_{1}^{4}$\textit{\ be a Cartan curve. Then }$\gamma $\textit{\ is a
pseudo-spherical curve iff }$k_{2}$\textit{\ is a nonzero constant.}
\end{Theorem}

\begin{Theorem}\textit{A Cartan curve }$\gamma \left(
t\right) $\textit{\ in }$%
\mathbb{R}
_{1}^{4}$\textit{\ fully lies on a pseudo-sphere iff there exists a fixed
point }$A$\textit{\ such that for each }$t\in I$\textit{, }$\left\langle
A-\gamma \left( t\right) ,\gamma ^{\prime }\left( t\right) \right\rangle =0.$
\end{Theorem}

\section{Bertrand Curve Corresponding to A Spacelike Geodesic on S$_{1}^{3}$}

\begin{Theorem}\textit{Let }$\mathit{\gamma }$ \textit{be a
spacelike geodesic curve on }$S_{1}^{3}.$\textit{\ Then,} \textit{%
\begin{equation*}
\tilde{\gamma}\left( s\right) =a\int \gamma \left( \upsilon \right)
d\upsilon +a\coth \theta \int e\left( \upsilon \right) d\upsilon +c
\end{equation*}%
is a Bertrand curve where }$\mathit{a}$\textit{\ and }$\theta $ \textit{are
constant numbers, }$\mathit{c}$\textit{\ is a constant vector.}
\end{Theorem}

\noindent \textbf{Proof. }We will use the frame $\{\gamma
(s),t(s),n(s),e(s)\}$ of $\gamma $ given in the previous section. In this
frame, let we choose $e\left( s\right) $ as a timelike vector (If $e\left(
s\right) $ is a spacelike vector, the proof is similar). Hence $n\left(
s\right) $ is spacelike and $\delta (\gamma (s))=-1.$ Using the equation
(1), we can easily calculate that%
\begin{eqnarray*}
\tilde{\gamma}^{\prime }\left( s\right)  &=&a\left[ \gamma \left( s\right)
+\coth \theta e\left( s\right) \right]  \\
\tilde{\gamma}^{\prime \prime }\left( s\right)  &=&a\left[ t\left( s\right)
+\coth \theta \tau _{g}\left( s\right) n\left( s\right) \right]  \\
\tilde{\gamma}^{\prime \prime \prime }\left( s\right)  &=&a[-\gamma \left(
s\right) +\delta (\gamma (s))\kappa _{g}\left( s\right) \tau _{g}\left(
s\right) t\left( s\right)  \\
&&+\left( \kappa _{g}\left( s\right) +\coth \theta \tau _{g}^{\prime }\left(
s\right) \right) n\left( s\right) +\coth \theta \tau _{g}^{2}\left( s\right)
e\left( s\right) ]
\end{eqnarray*}%
Since $\left\langle \tilde{\gamma}^{\prime }\left( s\right) ,\tilde{\gamma}%
^{\prime }\left( s\right) \right\rangle =-\dfrac{a%
{{}^2}%
}{\sinh
{{}^2}%
\theta }$, the curve $\tilde{\gamma}$ is timelike. If we calculate the first
and second curvatures of $\tilde{\gamma}$ by using the equations in $\left[ 8%
\right] $, we have%
\begin{eqnarray*}
\kappa \left( s\right)  &=&\frac{\sinh
{{}^2}%
\theta \sqrt{1+\coth
{{}^2}%
\theta \tau _{g}^{2}}}{a} \\
\tau \left( s\right)  &=&\frac{A\sinh \theta }{a\sqrt{1+\coth
{{}^2}%
\theta \tau _{g}^{2}}}
\end{eqnarray*}%
where $A=\sqrt{\cosh
{{}^2}%
\theta \left( \tau _{g}^{2}+1\right)
{{}^2}%
-\kappa _{g}^{2}\left( 1+\coth
{{}^2}%
\theta \tau _{g}^{2}\right) }.$ Since $\tau _{g}$ and $\kappa _{g}$ are
constants, we can choose $\beta =\dfrac{-a\sinh \theta \sqrt{1+\coth
{{}^2}%
\theta \tau _{g}^{2}}}{A}$ and $\alpha =\dfrac{a\coth
{{}^2}%
\theta }{\sqrt{1+\coth
{{}^2}%
\theta \tau _{g}^{2}}}$, then we have $\alpha \kappa +\beta \tau =1$. Hence $%
\tilde{\gamma}$ is a Bertrand curve.

\section{Bertrand Curve Corresponding to A Null Helix on S$_{1}^{3}$}

\begin{Theorem} \textit{Let }$\mathit{\gamma }$ \textit{be a null
helix on }$S_{1}^{3}.$\textit{\ Then,} \textit{%
\begin{equation*}
\tilde{\gamma}\left( s\right) =a\int L\left( \upsilon \right) d\upsilon
+a\coth \theta \int W_{2}\left( \upsilon \right) d\upsilon +c
\end{equation*}%
is a Bertrand curve where }$\mathit{a}$\textit{\ and }$\theta $ \textit{are
constant numbers, }$\mathit{c}$\textit{\ is a constant vector.}
\end{Theorem}

\noindent \textbf{Proof. }%
\begin{eqnarray*}
\tilde{\gamma}^{\prime }\left( t\right) &=&a\left[ L\left( s\right) +\coth
\theta W_{2}\left( t\right) \right] \\
\tilde{\gamma}^{\prime \prime }\left( t\right) &=&a\left[ 1-\coth \theta
k_{2}\right] W_{1}\left( t\right) \\
\tilde{\gamma}^{\prime \prime \prime }\left( t\right) &=&a\left[ k_{1}\left(
\coth \theta -1\right) L\left( t\right) -\left( 1-\coth \theta k_{2}\right)
N\left( t\right) \right]
\end{eqnarray*}%
Since $\left\langle \tilde{\gamma}^{\prime }\left( t\right) ,\tilde{\gamma}%
^{\prime }\left( t\right) \right\rangle =a%
{{}^2}%
\coth
{{}^2}%
\theta $, the curve $\tilde{\gamma}$ is spacelike. If we calculate the first
and second curvatures of $\tilde{\gamma}$, we have%
\begin{eqnarray*}
\kappa \left( t\right) &=&\frac{\left( 1-\coth \theta k_{2}\right) }{a\coth
{{}^2}%
\theta } \\
\tau \left( t\right) &=&\frac{\sqrt{k_{1}^{2}\cosh
{{}^2}%
\theta -1}}{\cosh \theta }
\end{eqnarray*}
Since $k_{1}$ and $k_{2}$ are constants, we can choose $\beta =-\dfrac{\cosh
{{}^3}%
\theta }{\sqrt{k_{1}^{2}\cosh
{{}^2}%
\theta -1}}$ and $\alpha =\dfrac{a\cosh
{{}^2}%
\theta }{\left( 1-\coth \theta k_{2}\right) }$, then we have $\alpha \kappa
+\beta \tau =1$. Hence $\tilde{\gamma}$ is a Bertrand curve.

\textbf{Received: Month xx, 200x}

\end{document}